\documentclass[12pt]{article}

\usepackage{diagbox}
\usepackage{mathtools}
\usepackage{bbm}
\usepackage{latexsym}
\usepackage{epsfig}
\usepackage{amsmath,amsthm,amssymb,enumerate}

\usepackage[a-1b]{pdfx}
\usepackage{hyperref}

\usepackage{color}

\def\reals{\mathbb R}

\def\a{\alpha}  \def\d{\delta} 
\def\e{\varepsilon}    
  
\def\z{\zeta}     \def\l{\lambda}
  \def\n{\nu} \def\p{\pi}
\def\r{\rho}  \def\s{\sigma} 
\def\t{\tau}

\newtheorem{theorem}{Theorem}

\newcommand{\brac}[1]{\left(#1\right)}

\newcommand{\bfrac}[2]{\left(\frac{#1}{#2}\right)}

\def\cE{{\cal E}}

\newcommand{\set}[1]{\left\{#1\right\}}

\def\E{\mathbb{E}}
\def\Var{\mathbb{V}\mathbb{A}\mathbb{R}}
\def\Pr{\mathbb{P}}

\newcommand{\ignore}[1]{}

\def\cE{{\mathcal E}}

\def\cX{{\mathcal X}}
\def\cY{{\mathcal Y}}

\newcommand{\beq}[2]{\begin{equation}\label{#1}#2\end{equation}}

\usepackage{tikz}
\usetikzlibrary{decorations.pathmorphing}
\usetikzlibrary{positioning}
\usetikzlibrary{arrows,automata}
\usetikzlibrary{shapes.misc}
\usetikzlibrary{backgrounds}
\usetikzlibrary{arrows,shapes}

\begin{document}

\date{}
\title{Aspects of a randomly growing cluster in $\reals^d,d\geq 2$}
\author{Alan Frieze\thanks{Department of Mathematical Sciences, Carnegie Mellon University, Pittsburgh PA15213, USA. Research supported in part by NSF grant DMS1952285 Email: frieze@cmu.edu} \and Ravi Kannan\thanks{Simons Institute for the Theory of Computing, UC Berkeley, Berkeley CA94720-2190. Email: kannan100@gmail.com} \and Wesley Pegden\thanks{Department of Mathematical Sciences, Carnegie Mellon University, Pittsburgh PA15213, USA. Research supported in part by NSF grant DMS1700365.Email: wes@math.cmu.edu}}
\maketitle
\begin{abstract}
We consider a simple model of a growing cluster of points in $\reals^d,d\geq 2$. Beginning with a point $X_1$ located at the origin, we generate a random sequence of points $X_1,X_2,\ldots,X_i,\ldots,$. To generate $X_{i},i\geq 2$ we choose a uniform integer $j$ in $[i-1]=\set{1,2,\ldots,i-1}$ and then let $X_{i}=X_j+D_i$ where $D_i=(\d_1,\ldots,\d_d)$. Here the $\d_j$ are independent copies of the Normal distribution $N(0,\s_i)$, where $\s_i=i^{-\a}$ for some $\a>0$. We prove that for any $\alpha>0$ the resulting point set is bounded a.s., and moreover, that the points generated look like samples from a $\beta$-dimensional subset of $\reals^d$ from the standpoint of the minimum lengths of combinatorial structures on the point-sets, where $\beta=\min(d,1/\alpha)$.
\end{abstract}
\section{Introduction}
In this short note, we study the following process: beginning with a point $X_1$ located at the origin, we generate a random sequence of points $X_1,X_2,\ldots,X_i,\ldots,$ in $\reals^d$. To generate $X_{i},i\geq 2$ we choose a uniform integer $j$ in $[i-1]=\set{1,2,\ldots,i-1}$ and then let $X_{i}=X_j+D_i$ where $D_i=(\d_1,\ldots,\d_d)$. Here the $\d_j$ are independent copies of the Normal distribution $N(0,\s_i)$, where $\s_i=i^{-\a}$ for some $\a>0$.  Thus as more and more points are added, new points are likely to cluster around old points.

We denote the set points $\set{X_1,X_2,\ldots,X_n}$ by $\cX_n$ and by $\cX_\infty=\bigcup_{n=1}^\infty \cX_n$ the set of all points generated by the process.  
Our first result shows that there is an exponential tail on the diameter $\r(\cX)=\max\set{|X|:X\in \cX}$ of the resulting infinite cluster:
\begin{theorem}\label{th1}
$\Pr(\r(\cX_\infty)\geq L)\leq e^{-L^2/600d}$ for large $L$.
\end{theorem}
As a consequence, the convex hull of $\cX_\infty$ is bounded a.s., by Borel-Cantelli applied to the events $\{\r(\cX_\infty)\geq L\}$ for $L=1,2,\dots.$  This stands in contrast to the case of a the set $\cY_\infty=\{Y_1,Y_2,\dots\}$ where each $Y_i$ is an independent standard Gaussian in $\reals^d$, which is unbounded a.s.~(and everywhere dense).

Our next theorem concerns the length $L_n$ of the minimum spanning tree on this collection of points under the Euclidean distance.  We note that the length of such minimum Euclidean structures are not just relevant from an optimization standpoint but can also be seen as a way of capturing the dimensionality of a set or distribution:  For $n$ points chosen uniformly from a compact subset $\Omega\subseteq \reals^d$ of dimension $\beta$ (e.g., a $\beta$-dimensional manifold for integer $\beta$, or a suitably regular fractal of dimension $\beta$ for non-integer $\beta$), the expected length of a spanning tree through the points is grows like $n^{1-1/\beta}$  \cite{bright}.  We show that from the standpoint of the length of a minimum spanning tree, the points generated by the process we study here look like uniform samples from a subset of $\reals^d$ of dimension $\min(d,\frac 1 \alpha)$. $a_1,\ldots,a_6$ are absolute constants.
\begin{theorem}\label{th2}
\[
\L_n\text{ satisfies }\begin{cases}a_1n^{1-1/d}\leq \E(L_n)\leq a_2n^{1-1/d}&\a<1/d.\\a_3n^{1-\a}\leq \E(L_n)\leq a_4n^{1-\a}\log^3n&\a> 1/d.\\a_5(n/\log n)^{1-1/d}\leq \E(L_n)\leq a_6n^{1-1/d}&\a=1/d.
\end{cases}
\]
\end{theorem}
Note that the particular choice of \emph{spanning tree} as our combinatorial structure is not so important here.  Indeed if $T$ and $H$ are the lengths of the minimum spanning tree and Hamilton cycle on the point-set, respectively, then we have $T\leq H\leq 2T$ and so the statement of Theorem \ref{th2} holds immediately for Hamilton cycles in place of trees here.  For other spanning structues like 2-factors or perfect matchings (for even $n$), the upper bounds in the theorem translate immediately, and the proofs of our lower bounds translate as well; in particular our proofs show not just that the lower bounds in Theorem \ref{th2} apply to the length of a spanning tree on the points in $\cX_n$, but to the total length of any collection of edges of linear size among the points $\cX_n$.

\section{Maximum distance: proof of Theorem \ref{th1}}
We define a tree $T_n$ with vertex set $\cX_n$ and edges of the form $X_iX_{\p(i)}$ for $i\in [n]$. Thus if $X_i$ chooses to be ``close'' to $X_j$ then we add the edge $X_iX_j$ to to $T_n$.

It is important to note that $T_n$ has the structure of a {\em random recursive tree}, see for example Chapter 14.2 of Frieze and Karo\'nski \cite{FK}.

Let $\l(i)=\l_n(i)$ denote the level of $X_i$ in the tree $T_n$, i.e. the number of edges from $X_i$ to $X_1$ in $T_n$. Let $\cE(i,\ell,L)$ be the event that $\l(i)=\ell$ and that the length of the edge from $X_i$ to its parent in $T_n$ is at least $\tfrac{L}{\ell^2\z(2)}$, where $\z(2)=\sum_{k-1}^\infty k^{-2}=\p^2/6$. If none of these events occur then every $i$ is at distance at most $\sum_{\ell=1}^\infty\tfrac{L}{\ell^2\z(2)}=L$ from the origin $X_1$.

In general when $i\leq  m$ we have that for integers $t\leq m$,
\[
\Pr(\l_m(i)>t)\leq \sum_{\substack{S\subseteq [m]\\|S|=t}}\prod_{j\in S}\frac{1}{j}\leq \frac{1}{t!}\brac{\sum_{j=1}^m\frac{1}{j}}^t\leq \bfrac{e(1+\log m)}{t}^t.
\]
It follows that
\beq{0}{
\Pr(\l(i)\geq 10(1+\log m))\leq m^{-4},\quad \text{ for $m$ large.}
}
Now we have the following inequality for $N(0,\s)$:
\beq{Tail}{
\Pr(N(0,\s)\geq x)\leq \frac{\s e^{-x^2/2\s^2}}{x(2\p)^{1/2}}.
}
We see from \eqref{Tail} that
\beq{chest}{
\Pr(\cE(i,\ell,L))\leq d\Pr\brac{N(0,i^{-\a})\geq \frac{L}{d^{1/2}\ell^2\z(2)}}\leq \frac{d\ell^2\z(2)}{(2\p)^{1/2}Li^\a} \exp\set{-\frac{L^2i^{2\a}}{2d\ell^4\z(2)^2}}.
}
(If $(\d_1^2+\d_2^2+\cdots+\d_d^2)^{1/2}\geq u=L/(\ell^2\z(2))$ then there exists $i$ such that $\d_i\geq u/d^{1/2}$. We can make a small improvement by using the bound on the upper tail of the  $\chi^2$-distribution in Laurent and Massart \cite{BM}.)

So for $L\leq k_1<k_2\leq n$ we have, using \eqref{0} and \eqref{chest}, 
\begin{align*}
\Pr\brac{\exists i\in [k_1,k_2]:\;\bigcup_{\ell\leq 10\log i}\cE(i,\ell,L)\text{ occurs}}&\leq \sum_{i=k_1}^{k_2}\brac{i^{-4}+\sum_{\ell=1}^{10(1+\log i)}\exp\set{-\frac{L^2i^{2\a}}{3d\ell^4\z(2)^2}}}\\
&\leq (k_2-k_1)\exp\set{-\frac{L^2k_1^{2\a}}{4d(10\log k_2)^4\z(2)^2}}+k_1^{-3}.
\end{align*}

So, let $m_0=n$ and $m_t=\log^{4/\a}m_{t-1}$ for $t=1,2,\ldots,t_0=\min\set{t:m_t\leq M}$ where $M=e^{L^2/1000d}$. It follows that
\begin{align*}
\sum_{t=1}^{t_0}\Pr(\exists i\in [m_t,m_{t-1}]:\;\cE(i,\ell,L)\text{ occurs for some $\ell$})&\leq \sum_{t=1}^{t_0}\brac{(m_{t-1}-m_t) \exp\set{-\frac{L^2m_t^{2\a}}{4d(10m_t^{\a/4})^4\z(2)^2}}+m_t^{-3}}\\
&\leq 2\sum_{t=1}^{t_0}m_t^{-3}\leq \frac{1}{M^2}.
\end{align*}
It follows that 
\[
\Pr(\exists i:dist(i)\geq L)\leq \frac{1}{M^2}+\frac{1}{M^3}+ e^{-L^2/500d}\leq e^{-L^2/600d}.
\]
The term 
\[
\frac{1}{M^{3}}+\frac{Md(10(1+\log M))^2\z(2)e^{-L^2/300d}}{(2\p)^{1/2}L}\leq \frac{1}{M^3}+ e^{-L^2/500d}
\]
 arises from applying \eqref{Tail} (with $\s=1$) and \eqref{0} to bound the probability that $\cE(i,\ell,L)$ occurs for some $i,\ell\leq M$. This completes the proof of Theorem \ref{th1}.
\section{Minimum spanning tree}
\subsection{Upper bound}
We bound the length of the recursive tree $T_n$. Very crudely, the cost of the first $\log n$ edges is $O(\log^2n)$ q.s.\footnote{A sequence of events $\cE_n,n\geq 1$ occurs {\em quite surely} (q.s.) if $\Pr(\neg\cE_n)=O(n^{-K})$ for any constant $K>0$.} Next let $L_i=i^{-\a}\log^3n$. Suppose that $\cE(i,\ell,L_i)$ does not occur for $i\geq \log n$. Then the length of the tree produced is at most
\[
O(\log^2n)+\sum_{i=\log n}^n\sum_{\ell=1}^{10(1+\log i)}\frac{L_i}{\ell^2\z(2)}\leq O(\log^2n)+\sum_{i=\log n}^nL_i.
\]
We see from \eqref{chest} that the probability we fail to produce a tree of the claimed size is at most
\[
o(1)+\sum_{i=\log n}^n\sum_{\ell=1}^{10(1+\log i)}\frac{d\ell^2\z(2)}{\log^3n}\exp\set{-\frac{\log^6n}{d\ell^4\z(2)^2}}=o(1).
\]
Thus w.h.p. there is a tree of length at most
\[
O(\log^2n)+\sum_{i=\log n}^n\frac{\log^3n}{i^\a}\leq n^{1-\a}\log^3n.
\]
This gives the upper bound in Theorem \ref{th2} for $\a\geq 1/d$. For $\a<1/d$ we appeal to the fact the claimed upper bound holds for all sets of $n$ points, in a bounded region, see for example Steele and Snyder \cite{SS}. So from Theorem \ref{th1} we can claim that the expected length of the minimum spanning tree is at most
\[
c_3n^{(d-1)/d}\int_{L=0}^{\infty}Le^{-L^2/600d}dL=O(n^{(d-1)/d}). 
\]
This proves the upper bound for $\a<1/d$.
\subsection{Lower bounds}
Consider two vertices $i,j$ whose common ancestor in the recursive tree is $m$. Then we have
\[
\Pr(|X_i-X_j|\leq \d)\leq \Pr(|N(0,m^{-\a})-N(0,m^{-\a})|\leq \d)^d=\Pr(|N(0,2m^{-\a})|\leq \d)^d=O((\d m^{\a})^d).
\]
Now in general, the expected number of pairs $i,j$ with common ancestor $m$ is at most $2n^2/m^2$, see equation \eqref{nm} in Section \ref{urn}. So, if $Z_\d$ denotes the number of pairs of vertices at distance at most $\d$, then for some constants $C_1,C_2$, 
\[
\E(Z_\d)\leq C_1n^2\d^{d}\sum_{m=1}^nm^{\a d-2}\leq C_2n^2\d^{d}\times\begin{cases}1&\a<1/d\\n^{\a d-1}&\a>1/d\\ \log n&\a=1/d.\end{cases}
\]
If $\a<1/d$ then we can put $\d=\e n^{-1/d}$ for small $\e>0$ and see that the expected number of pairs $i,j$ at distance at most $\d$ is at most $C_2\e^{d}n$. In which case the expected length of the minimum spanning tree is at least
\beq{add1}{
((n-1)-C_2\e^{d}n)\d\geq c_1 n^{1-1/d}
}
for constant $c_1$.

If $\a>1/d$ then we can put $\d=\e n^{-\a}$ and see that the expected number of pairs $i,j$ at distance at most $\d$ is at most $C_2\e^{d}n$. In which case the expected length of the minimum spanning tree is at least
\beq{add0}{
((n-1)-C_2\e^{d}n)\d \geq c_2  n^{1-\a}
}
for constant $c_2$. 

If $\a=1/d$ we put $\d=\e(n\log n)^{-1/d}$ and see that the expected number of pairs $i,j$ at distance at most $\d$ is at most $2C_2\e^{d}n$. In which case the expected length of the minimum spanning tree is at least
\beq{add}{
((n-1)-2C_2\e^{d}n)\d \geq c_3  (n/\log n)^{1-1/d}
}
for constant $c_3$. 

This completes the proof of Theorem \ref{th2}. (Note that \eqref{add1}, \eqref{add0}, \eqref{add} show that the lower bounds in Theorem \ref{th2} apply to any set of $\Omega(n)$ edges.)

\subsubsection{Polya-Eggenburger Urn}\label{urn}
In the Polya-Eggenburger Urn with parameters $W_0=1,B_0=m-1,\t_0=W_0+B_0,s=1$ we start with an urn containing $W_0$ white balls, $B_0$ blue balls. At each round we choose a ball at random and replace it and then add $s$ balls of the same color to the urn. For us, the balls are the vertices of the tree. The white balls are the descendants of vertex $m$. Let $W_n$ denote the number of white balls in the urn after $n$ rounds. Then Corollary 5.1.1 of \cite{M} states
\[
\E(W_n)=\frac{W_0}{\t_0}sn+W_0\text{ and }\Var(W_n)=\frac{W_0B_0s^2n(sn+\t_0)}{\t_0^2(\t_0+s)}.
\]
Plugging in our values, we get $\E(W_{n-m})=(n-m)/m+1=n/m$ and 
\beq{nm}{
\E(W_{n-m}^2)=\Var(W_{n-m})+\E(W_{n-m})^2=\frac{(m-1)n(n-m)}{m^2(m+1)}+\bfrac{n}{m}^2\leq \frac{2n^2}{m^2}.
}
\section{Summary}
We have introduced a new model of a point process and have proved bounds on its spread and the cost of the minimum spanning tree through the points. We could have considered starting the process with $k>1$ points placed arbitrarily. This would involve $k$ trees with sizes determined by the Polya-Eggenburger model and it is not hard to see that our two theorems are still valid. It might be of some interest to try and remove the polylog factors from Theorem \ref{th2}. Maybe also, one could try other sequences of standard deviation, other than $i^{-\a}$.

One natural question to ask, is as to what happens when $\a=0$, i.e. when the $\d_j$ in the definition of the $X_i$ are $N(0,1)$. In this case Theorem \ref{th1} fails. We know that w.h.p. the depth of $T_n$ is $\Omega(\log n)$. In which case, the distance of leaves in $T_n$ from the root $X_1$ are bounded below by the sum of $\Omega(\log n)$ standard normals and so they will w.h.p. be $\Omega(\log n)$ from $X_1$.

\end{document}